\documentclass{amsart}
\usepackage{hyperref}
\usepackage[backrefs]{amsrefs}

\newtheorem{theorem}{Theorem}[section]

\newtheorem{algorithm}[theorem]{Algorithm}

\theoremstyle{definition}
\newtheorem{definition}[theorem]{Definition}

\theoremstyle{remark}

\numberwithin{equation}{section}

\newcommand{\Pro}{\mathbb{P}}

\newcommand{\Aut}{\operatorname{Aut}}
\newcommand{\Sym}{\operatorname{Sym}}
\newcommand{\GL}{\operatorname{GL}}
\newcommand{\Order}{\operatorname{Order}}
\newcommand{\Id}{\operatorname{Id}}
\newcommand{\rank}{\operatorname{rank}}

\begin{document}

\title{Equations of Riemann surfaces with automorphisms}


\author{David Swinarski}
\address{Department of Mathematics\\ 
Fordham University\\ 
113 W 60th St Room 813\\ 
New York, NY 10023}
\email{dswinarski@fordham.edu}


\date{}

\begin{abstract}
We present an algorithm for computing equations of canonically embedded
Riemann surfaces with automorphisms.    A variant of this
algorithm with many heuristic improvements is used to produce equations of Riemann surfaces $X$ with large
automorphism groups (that is, $|\Aut(X)| > 4(g_X-1)$) for
genus $4 \leq g_X \leq 7$.  The main tools are the Eichler trace
formula for the character of the action of $\Aut(X)$ on holomorphic
differentials, algorithms for producing matrix generators of a
representation of a finite group with a specified irreducible character, and Gr\"obner basis techniques for computing flattening
stratifications.  
\end{abstract}

\maketitle

Riemann surfaces (or algebraic curves) with automorphisms have been important objects of study
in complex analysis, algebraic geometry, number theory, and
theoretical physics for over a century, as their symmetries often
permit us to do calculations that would otherwise be intractable.  

Such Riemann surfaces are special in the sense that a general
Riemann surface of genus $g \geq 3$ has no nontrivial automorphisms.  Moreover,
the group of automorphisms of a Riemann surface of genus $g \geq 2$ is
finite.

Breuer and Conder performed computer searches that for each genus $g$ list the
Riemann surfaces of genus $g$ with large automorphism groups (that is,
$|\Aut(X)| > 4(g_X-1)$).
Specifically, they list sets of surface kernel generators
(see Definition \ref{surface kernel definition} below), which describe these Riemann surfaces as branched covers of
$\Pro^1$.  Breuer's list extends to genus $g=48$, and Conder's list
extends to genus $g=101$ \cites{Breuer,Conder}.   Even for small values of $g$, these lists are extremely large, as a surface $X$ may appear several times for various subgroups of its full automorphism group.  In \cite{MSSV}, Magaard, Shaska, Shpectorov, and V\"olklein refined Breuer's list by determining which surface kernel generators
correspond to the full automorphism group of the Riemann surface.

To my knowledge, at this time there is no general algorithm published in the literature
for producing equations of these Riemann surfaces under any embedding
from this data.  Here, I present an algorithm to compute 
canonical equations of nonhyperelliptic Riemann surfaces with automorphisms.  The main tools are the Eichler trace
formula for the character of the action of $\Aut(X)$ on holomorphic
differentials, algorithms for producing matrix generators of a
representation of a finite group with a specified irreducible
character, and Gr\"obner basis techniques for computing flattening
stratifications.  
A variant of this algorithm with many heuristic improvements is used
to produce equations of the nonhyperelliptic Riemann
surfaces with genus $4 \leq g_X \leq 7$ satisfying $|\Aut(X)| > 4(g_X-1)$.

Here is an outline of the paper.  In Section \ref{algorithm
  section}, I describe the main algorithm.  In Section
\ref{heuristics section}, I describe several heuristics that simplify or speed up the main algorithm.  In Section \ref{example
  section}, I describe one example in detail, a genus 7 Riemann
surface with 64 automorphisms.  In Section \ref{results
  section}, I give  equations of
selected canonically embedded Riemann surfaces with $4 \leq g_X \leq
7$ along with matrix surface kernel generators.  

\subsection*{Acknowledgements}
A large number of people have supplied advice and encouragement on
this project over a period of many years. The undergraduates from my Introductory
 VIGRE Research Group held in Fall 2010 at the University of Georgia
 and my student Darcy Chanin performed many calculations for 
 genus 4, 5, and 6 surfaces.  I am grateful to the computational algebra
 group at the University of Sydney, especially John Cannon and Mark Watkins, 
for hosting me for a visit in June 2011 where I began programming the
main algorithm in \texttt{Magma} \cite{Magma}.   I have had many helpful
conversations with my classmates and colleagues at Columbia
University, the University of Georgia, and Fordham University.  Valery Alexeev and James McKernan suggested the
algorithm for matrix generators of representations outlined in Section
\ref{matrix generators section}.  
Finally, I am grateful to Jennifer Paulhus, Tony Shaska, and
John Voight, whose encouragement was essential in completing this project.  This work was
partially supported by the University of Georgia's NSF VIGRE grant
DMS-03040000, a Simons Foundation Travel Grant, and a Fordham Faculty Research Grant.

\subsection*{Online material}
My webpage for this project is \cite{mywebpage}.  This page contains links to the latest version of my \texttt{Magma} 
code, files detailing the calculations for specific examples, and many equations that are omitted in the tables in Section \ref{results
  section}.

In future work, Jennifer Paulhus and I plan to include much of the data
described in this paper and on the website \cite{mywebpage}  in
the L-Functions and Modular Forms Database at \texttt{lmfdb.org}.

\section{The main algorithm} \label{algorithm section}

We begin by stating the main algorithm.  Then, in the following subsections, we
discuss each step in more detail, including precise definitions
and references for terms and facts that are not commonly known.  

\begin{algorithm} \label{main algorithm} \mbox{} \\
\textsc{Inputs:} 
\begin{enumerate}
\item A finite group $G$; 
\item an integer $g \geq 2$; 
\item a set
of surface kernel generators $(a_1,\ldots,a_{g_0}; b_1,\ldots,b_{g_0};
g_1,\ldots,g_r)$ determining a family
of nonhyperelliptic Riemann
surfaces $X$ of genus $g$ with $G \subset \Aut(X)$\\
\end{enumerate}
\textsc{Output:} A locally closed set $B \subset \mathbb{A}^{n}$ and
a family of smooth curves $\mathcal{X} \subset \Pro^{g-1} \times B$ such that for
each closed point $b \in B$, the fiber $\mathcal{X}_b$ is a smooth
genus $g$ canonically embedded curve with $G \subset \Aut(\mathcal{X}_b)$.\\

\begin{enumerate}
\item[Step 1.] Compute the conjugacy classes and character table of $G$.
\item[Step 2.] Use the Eichler trace formula to compute the character of the
  action on differentials and on cubics in the canonical
  ideal.  
\item[Step 3.] Obtain matrix generators for the action on holomorphic differentials.
\item[Step 4.] Use the projection formula to obtain candidate cubics.
\item[Step 5.] Compute a flattening stratification and select the locus yielding 
  smooth algebraic curves with degree $2g-2$ and genus $g$.  
\end{enumerate}
\end{algorithm}

\subsection{Step 1: conjugacy classes and character table of $G$}

This step is purely for bookkeeping.  It is customary to list the
conjugacy classes of $G$ in increasing order, and to list the rows in
a character table by increasing degree.  However, there is no
canonical order to either the conjugacy classes or the irreducible characters.  Given two different descriptions of a finite group $G$, modern
software such as \texttt{Magma} may order the classes or the
irreducible characters of $G$ differently.  Hence, we compute and fix
these at the beginning of the calculation.

\subsection{Step 2: Counting fixed points and the
  Eichler trace formula}
Here we define surface kernel generators for the automorphism
group of a Riemann surface.   These generators determine the Riemann
surface as a branched cover of $\Pro^1$ and are used in a key formula
(see Theorem \ref{count fixed points} below)
for counting the number of fixed points of an automorphism.

\begin{definition}[cf.~\cite{Breuer} Theorem 3.2, Theorem 3.14]  
\label{surface kernel definition}
A \emph{signature} is a list of integers $(g_0; e_1,\ldots,e_r)$ with
$g_0\geq 0$, $r \geq 0$, and $e_i \geq 2$.  \\

A set of \emph{surface kernel generators} for a finite group $G$ and signature
$(g_0; e_1,\ldots,e_r)$ is a sequence of elements $a_1,\ldots,
a_{g_0}, b_{1},\ldots, b_{g_0}, g_1,\ldots, g_r \in G$ such that
\begin{enumerate}
\item  $\langle a_1,\ldots,
a_{g_0}, b_{1},\ldots, b_{g_0}, g_1,\ldots, g_r \rangle = G$;
\item $\Order(g_i)=e_i$; and
\item $\prod_{j=1}^{g_0} [a_j,b_j] \prod_{i=1}^{r} g_i = \Id_G$.
\end{enumerate}
\end{definition}

Surface kernel generators have many other names in other papers; they are
called \emph{ramification types} in \cite{MSSV} and \emph{generating
  vectors} in \cite{Paulhus}.  

As explained in \cite{Breuer}*{Section 3.11}, surface kernel
generators describe the quotient morphism $X \rightarrow X/G$ as a
branched cover.  Here $X$ is a Riemann surface of genus $g$, $G$ is a
subgroup of $\Aut(X)$, the quotient $X/G$ has genus $g_0$, the quotient
morphism branches over $r$ points, and the integers $e_i$
describe the ramification over the branch points.

In the sequel we will be primarily interested in large automorphism
groups, that is, $|\Aut(X)| > 4(g_X-1)$.  In this case, the Riemann-Hurwitz
formula implies that $g_0 = 0$ and $3 \leq r \leq 4$.

Surface kernel generators are used in the following formula for 
the number of fixed points of an automorphism:

\begin{theorem}[\cite{Breuer}*{Lemma 11.5}] \label{count
    fixed points}
Let $\sigma$ be an automorphism of order $h>1$ of a Riemann surface $X$ of genus
$g\geq 2$.  Let $(g_1,\ldots,g_r)$ be part of a set of surface kernel
generators for $X$, and let $(m_1,\ldots, m_r)$ be the orders of these elements.  Let $\mbox{Fix}_{X,u}(\sigma)$ be the set of fixed points of
$X$ where $\sigma$ acts on a neighborhood of the fixed point by $z \mapsto
\exp(2 \pi i u/h) z$.  Then

\begin{displaymath}
|\operatorname{Fix}_{X,u}(\sigma)| = |C_{G}(\sigma)|
\sum_{\substack{g_i \, s.t. \\ h|m_i \\ \sigma \sim g_i^{m_i u/h} }}  \frac{1}{m_i}
\end{displaymath}

Here $C_{G}(\sigma)$ is the centralizer of $\sigma$ in $G$, and $\sim$
denotes conjugacy.
\end{theorem}

Next we recall the Eichler Trace Formula.  For a Riemann surface $X$,
let $\Omega_X$ be the holomorphic cotangent bundle, and let $\omega_X
= \bigwedge \Omega_X$ be the sheaf of holomorphic differentials.  The Eichler Trace Formula gives the character of
the action of $\Aut(X)$ on $\Gamma(\omega_X^{\otimes d})$.  

\begin{theorem}[Eichler Trace Formula
  \cite{FarkasKra}*{Theorem V.2.9}]  \label{Eichler Trace Formula} Suppose $g_{X} \geq 2$, and let $\sigma$ be a
  nontrivial automorphism of $X$ of order $h$.  Write $\chi_{d}$ for the character of the
   representation of $\Aut(X)$ on $\Gamma(\omega_X^{\otimes d})$.  Then
\begin{displaymath}
  \chi_d(\sigma) = \left\{
\begin{array}{ll}
1 + \displaystyle \sum_{\substack{1 \leq u < h \\ (u,h)=1}} |\operatorname{Fix}_{X,u}(\sigma)
  | \frac{\zeta_{h}^{u}}{1-\zeta_{h}^{u}}  & \mbox{ if $d=1$} \\
\displaystyle \sum_{\substack{1 \leq u < h \\ (u,h)=1}} |\operatorname{Fix}_{X,u}(\sigma)
  | \frac{\zeta_{h}^{u (d\%h)}}{1-\zeta_{h}^{u}}  & \mbox{ if $d \geq 2$} 
\end{array} \right.
\end{displaymath}
\end{theorem}

Together, the previous two results give a group-theoretic method for
computing the character of the $\Aut(X)$ action on $\Gamma(\omega_X^{\otimes
  d})$ starting from a set of surface kernel generators.  

We can use the character of $\Aut(X)$ on $\Gamma(\omega_X^{\otimes
  d})$ to obtain the character of $\Aut(X)$ on quadric and cubics in
the canonical ideal as follows.  Let $S$ be the 
coordinate ring of $\Pro^{g-1}$, let $I \subset S$ be the canonical ideal, 
and let $S_d$ and $I_d$ denote the degree $d$ subspaces of $S$ and $I$.  

By Noether's Theorem, the sequence 
\[
0 \rightarrow I_d \rightarrow S_d \rightarrow \Gamma(\omega_X^{\otimes d}) \rightarrow 0
\]
is exact for each $d \geq 2$, and by Petri's Theorem, the canonical ideal is generated either by quadrics or by quadrics and cubics.    
Thus, beginning with the character of the action on $\Gamma(\omega_X) \cong S_1$, we may compute 
the characters of the actions on $S_2 = \Sym^2 S_1$ and $S_3 = \Sym^3 S_1$ and $\Gamma(\omega_X^{\otimes 2})$ and 
$\Gamma(\omega_X^{\otimes 3})$, and then obtain the characters of the actions on $I_2$ and $I_3$.

\subsection{Step 3: matrix generators for a specified irreducible
  character} \label{matrix generators section}
From Step 2 we have the character of the action on $\Gamma(\omega_X)$.  We seek matrix generators for this action.  
It suffices to find matrix generators for each irreducible $G$-module appearing in $\Gamma(\omega_X)$.  

Given a finite group $G$ and an irreducible character $\chi$ of $G$, 
software such as \texttt{GAP} \cite{GAP} and \texttt{Magma} contain
commands for producing matrix generators of a representation $V$ of
$G$ with character $\chi$.  Finding efficient algorithms to produce
matrix generators with good properties (for instance, sparse matrices,
or matrices whose entries have small height, or matrices whose entries belong
to a low degree extension of $\mathbb{Q}$) is a subject of ongoing research \cites{Dabbaghian,DabbaghianDixon}.  It seems
that computer algebra systems implement several different algorithms that cover
many special cases.

I do not know a reference for a general algorithm.  Hence, I briefly present an algorithm that was suggested to
me by Valery Alexeev and James McKernan.  
This algorithm is not expected to perform efficiently;  it is included
merely to 
establish that Step 3 in Algorithm \ref{main algorithm} can be performed
algorithmically.  

\begin{algorithm} \mbox{} \\
\textsc{Inputs:} 
\begin{enumerate}
\item a finite group $G$ with generators $g_1,\ldots,g_r$; 
\item an irreducible character $\chi: G
\rightarrow \mathbb{C}$ of degree $n$.
\end{enumerate}
\textsc{Output:} matrices $M_{1},\ldots,M_{r} \in
\operatorname{GL}(n,\mathbb{C})$ such that the homomorphism $g_i
\mapsto M_i$ is a representation with character $\chi$
\begin{enumerate}
\item[Step 1.] Compute matrix generators for the regular representation $V$ of
  $G$.  These matrices are permutation matrices, and hence their entries are in $\{0,1\}$.  
\item[Step 2.] Use the projection formula (see Theorem \ref{projection formula
    theorem} below) to compute matrix generators $\rho_W(g)$ for a
  representation $W$ with character $n\chi$.  Let $K $ be
  the smallest field containing $\{ \chi(g) : g \in
  G\}$.  Note that $\mathbb{Q} \subseteq K \subseteq
  \mathbb{Q}[\zeta_{|G|}]$.  Then the matrix generators $\rho_W(g)$ lie in $\GL(n^2,K)$.
\item[Step 3.] Let $x_1,\ldots,x_{n^2}$ be indeterminates.  Let $M$ be the $|G|
  \times n^2$ matrix over $K$ whose rows are given by the vectors
  $\rho_W(g).(x_1,\ldots,x_{n^2})$. Let $X \subset \Pro^{n^2-1}_{K}$ be the determinantal variety $\rank
  M \leq n$.  Since representations of finite groups are completely
  reducible in characteristic zero, the 
  representation $W$ is isomorphic over $K$ to the direct sum
  $V_{\chi}^{\oplus n}$, and therefore $X(K)$ is non-empty.
\item[Step 4.] Intersect $X$ with generic hyperplanes with
  coefficients in $K$ to obtain a
  zero-dimensional variety $Y$.  
\item[Step 5.] If necessary, pass to a finite field extension $L$ of $K$ to obtain a reduced closed point
  $y \in Y(L)$.  
\item[Step 6.] The point $y$ (regarded as a vector in $W \otimes L$) generates the desired representation.
\end{enumerate}
\end{algorithm}

An example where this algorithm is used to produce matrix generators for the
degree two irreducible representation of the symmetric group $S_3$ is available at my webpage \cite{mywebpage}.  

Finally, we note that in \cite{Streit}, Streit describes a method for producing
matrix generators for the action of $\Aut(X)$ on $\Gamma(\omega_X)$
for some Bely\u{i} curves.  

\subsection{Step 4: the projection formula} 
Recall the projection formula for representations of finite groups.  
(See for instance \cite{FultonHarris} formula (2.31)).  

\begin{theorem}[Projection formula] \label{projection formula theorem}
Let $V$ be a finite-dimensional representation of a finite group $G$
over $\mathbb{C}$.
Let $V_1,\ldots,V_k$ be the irreducible representations of $G$, let
$\chi_i$ be their characters, and
let $V \cong \bigoplus_{i=1}^{k} V_{i}^{\oplus m_i}$.  Let $\pi: V \rightarrow V_{i}^{\oplus m_i}$
be the projection onto the $i^{th}$ isotypical component of $V$.  Then 
\[
\pi_i = \frac{\dim(V_i)}{|G|} \sum_{g \in G} \overline{\chi_{i}(g)} g.
\]
\end{theorem}
From Step 3, we have matrix generators for the $G$ action on $\Gamma(\omega_X) = S_1$.  Thus, we can compute 
matrix generators for the actions on $S_2$ and $S_3$, and use the projection formula to compute 
the isotypical subspace $S_{d,p}$ of degree $d$ polynomials on which $G$ acts with character $\chi_p$.  
In some a few examples, we have $I_{d,p} = S_{d,p}$, but more commonly, we have strict containment 
 $I_{d,p} \subset S_{d,p}$.  In this case we write elements of  $I_{d,p} $ as generic linear combinations of the basis elements 
of $S_{d,p}$ and then seek coefficients that yield a smooth algebraic curve with the correct degree and genus.  

The coefficients used to form these generic linear combinations form
the base space $\mathbb{A}^n$ of the family $\mathcal{X}$ produced by
the main algorithm.  

\subsection{Step 5: Flattening
  stratifications} \label{flattening stratification section}

\begin{theorem}[\cite{MumfordCurves}*{Lecture 8}]  Let $f: X \rightarrow S$
  be a projective morphism 
  with $S$ a reduced Noetherian scheme.  Then there exist locally closed subsets $S_1,\ldots, S_n$ such that $S =
  \sqcup_{i=1}^{n} S_i$ and $f|_{f^{-1}(S_i)} $ is flat.  
\end{theorem}

The stratification $S = \sqcup_{i=1}^{n} S_i$ is called a flattening
stratification for the map $f$.  Since $S$ is reduced, flatness
implies that over each stratum, the Hilbert polynomial of the fibers
is constant.  We find the stratum with Hilbert polynomial $P(t)
= (2g-2)t -g+1$, then intersect this stratum with the locus where the fibers
are smooth.  This completes the algorithm.

Flattening stratifications have been an important tool in theoretical algebraic
geometry for over 50 years.  There exist Gr\"obner
basis techniques for computing flattening stratifications; in the
computational literature, these are typically called
\emph{comprehensive} or \emph{parametric Gr\"obner bases}, or \emph{Gr\"obner systems}.
The foundational work on this problem was begun by Weispfenning, and
many authors, including Manubens and Montes, Suzuki and Sato, Nabeshima,
and Kapur, Sun, and Wang, have made important improvements on the original
algorithm \cites{Weispfenning,Nabeshima,KapurSunWang}.

The size of a Gr\"obner basis can grow very quickly with the number of
variables and generators of an ideal, and unfortunately, even the most
recent software cannot compute  flattening
stratifications for the examples we consider.  Thus, in section
\ref{partial flattening stratification section} below, we discuss a strategy
for circumventing this obstacle.

\section{Heuristic improvements} \label{heuristics section} 

Many steps of Algorithm \ref{main algorithm} can be run using a
computer algebra system, but even for modest examples, 
the flattening stratification required in the final step is intractable.  Therefore we discuss various 
heuristics that can be employed to speed the computation.  

\subsection{Tests for gonality and reduction to quadrics}
Given a set of surface kernel generators, it is useful to discover as early as possible whether the 
corresponding Riemann surface is hyperelliptic, trigonal, a plane quintic, or none of these.  We discuss these properties 
in turn.   

\textit{Hyperelliptic Riemann surfaces}.  Algorithm \ref{main algorithm}  supposes that one begins with surface kernel generators corresponding to a 
nonhyperelliptic curve.  However, we can easily test for
hyperellipticity if this  property  is not known in advance.  A Riemann surface $X$ is
  hyperelliptic if and only if $\Aut(X)$ contains a central involution
  with $2g_X+2$ fixed points. Thus, given a set of surface kernel generators, we can search for a central 
involution and count its fixed points using Theorem \ref{count fixed
  points} (or even better, using \cite{Breuer}*{Lemma 10.4}).  

In \cite{Shaska}, Shaska gives equations of the form $y^2=f(x)$ for hyperelliptic curves with automorphisms.  
Additionally, we can use the algorithm described in
\cite{StevensDeformations} to get the equations of $C$ under a linear
series such as the transcanonical embedding or bicanonical embedding.   

So suppose the Riemann surface is not hyperelliptic.  By Petri's Theorem, 
the canonical ideal is generated by quadrics if $X$ is not hyperelliptic, not trigonal, and not a plane quintic.  
Thus, ruling out these possibilities allows us to work with quadrics
instead of cubics, which significantly speeds up the algorithm.  This
leads us to consider trigonal Riemann surfaces and plane quintics.

\textit{Trigonal Riemann surfaces.}  
Trigonal Riemann surfaces may be divided into two types: cyclic trigonal and general trigonal \cite{CostaIzquierdo}.  
Cyclic trigonal curves can be detected by searching for degree three elements fixing $g+2$ points.  
Their automorphism groups have been classified \cite{BCGCyclicTrigonal}*{Theorem 2.1}, and one may hope for a paper treating
equations of cyclic trigonal Riemann surfaces as the paper 
\cite{Shaska} treats equations of hyperelliptic Riemann surfaces.

Less is known about general trigonal curves.  We have Arakawa's bounds
\cite{Arakawa}*{Remark 5}
and a few additional necessary conditions 
\cite{CostaIzquierdo}*{Prop. 4 and Lemma 5}.  We will not say more
about general trigonal Riemann surfaces here because after studying the
Riemann surfaces with large automorphism groups with genus $ 4 \leq g
\leq 7$, we learn \textit{a posteriori} that very few of them are general trigonal.

\textit{Plane quintics.}  Plane quintics only occur in genus 6, and the canonical model of a
plane quintic lies on the Veronese surface in $\Pro^5$.  Thus, we have a necessary condition: 
$X$ is a plane quintic only if $\Gamma(\omega_X) \cong \Sym^2 V$ for some (possibly reducible)
three-dimensional representation $V$ of $G$.  In practice, it is generally quite fast to discover whether a nonhyperelliptic 
non-cyclic trigonal genus 6 Riemann surface is a plane quintic.  

\subsection{Partial flattening stratifications} \label{partial
  flattening stratification section}
In this section we use several notions from the theory of Gr\"obner
bases.  We will not recall all the definitions here, and instead refer
to \cite{Eisenbud}*{Chapter 15} for the details.  

The algorithms for comprehensive Gr\"obner bases described in Section 
\ref{flattening stratification section} all begin with the same observation.  
Let $S=K[x_0,\ldots,x_m]$ be a polynomial ring over a field.  Let
$\preceq $ be a multiplicative term order on $S$.  Then a theorem of
Macaulay states that the Hilbert
function of $I$ is the same as the Hilbert function of its initial
ideal with respect to this term order
(see \cite{Eisenbud}*{Theorem 15.26}).  

Therefore, whenever two ideals in $S$ have Gr\"obner bases with the same
leading monomials with respect to some term order, they will have the
same initial ideal  for that term order, hence they must have the
same Hilbert function and Hilbert polynomial, and therefore they will
lie in the same stratum of a flattening stratification.  To reach a different 
stratum of the flattening stratification, it is necessary to alter the leading 
terms of the Gr\"obner basis --- for instance, by restricting to the 
locus where that coefficient vanishes.  

Here is a brief example to illustrate this idea.  Let $\mathbb{A}^2$
have coordinates $c_1,c_2$, and let $\Pro^3$ have coordinates
$x_0$,$x_1$,$x_2$,$x_3$. The ideal
\[
I = \langle c_1 x_0x_2-c_2 x_1^2, c_1 x_0x_3-c_2 x_1x_2, c_1
x_1x_3-c_2 x_2^2 \rangle
\]
defines a 2-parameter family of subschemes of $\Pro^3$.  A Gr\"obner basis for $I$ in
$\mathbb{C}[c_1,c_2][x_0,x_1,x_2,x_3]$ with respect to the lexicographic term order is
\[
\begin{array}{cc}
c_1 x_0 x_2-c_2 x_1^2, c_1 x_0 x_3-c_2 x_1 x_2, c_1 x_1 x_3-c_2 x_2^2,\\
(c_1 c_2-c_2^2) x_1 x_2^2, c_2 x_1^2 x_3-c_2 x_1 x_2^2, (c_1 c_2-c_2^2) x_1^2 x_2, c_2 x_0 x_2^2-c_2 x_1^2 x_2, \\
(c_1 c_2^2-c_2^3) x_1^4, (c_1 c_2^2-c_2^3) x_2^4, c_2^2 x_1 x_2^2 x_3-c_2^2 x_2^4, c_2^2 x_0 x_1^2 x_2-c_2^2 x_1^4.
\end{array}
\] 
Over the locus where $c_1$, $c_2$, and $c_1-c_2$ are invertible, the initial ideal
  is  \\$\langle x_0 x_2, x_0 x_3, x_1 x_3,x_1 x_2^2,x_1^2 x_2,x_1^4,x_2^4\rangle$ with
  Hilbert polynomial $P(t) = 8$.  On the other hand, when $c_1=0$, or $c_2=0$, or $c_1-c_2=0$, we  get a different initial ideal and Hilbert polynomial.  For example, the locus $c_1=c_2 \neq 0$ yields the twisted cubic with $P(t) = 3t+1$.

Note that to discover this locus, it is not necessary to compute the entire Gr\"obner basis; it would suffice for instance to compute the S-pair reduction for the first two generators, which yields $(c_1 c_2-c_2^2) x_1 x_2^2$.  

Modern software packages by Nabeshima, Montes, and Kapur, Sun, and
Wang can completely analyze this example.  However, these packages did not yield answers on the problems that arose in this work.  Therefore, I used the strategy outlined above. I partially computed a Gr\"obner basis in \texttt{Macaulay2} \cite{Macaulay2}, and set some coefficients to zero.  Remarkably, this was sufficient to obtain the equations of the genus $4 \leq g \leq 7$ Riemann surfaces with large automorphism groups.  Some of the families analyzed in this manner had as many as six coefficients $c_1,\ldots,c_6$.

\section{Example: a genus 7 Riemann surface with 64 automorphisms} \label{example section}
Magaard, Shaska, Shpectorov, and V\"olklein's tables show that there exists a
smooth, compact genus 7 Riemann surface with
automorphism group $G$ given by the group labeled $(64,41)$ in the GAP
library of small finite groups.  It has $X/G \cong \mathbb{P}^1$.  The quotient
morphism is branched over 3 points of $\mathbb{P}^1$, and the
ramification indices over these points are  2, 4, and 16.  

A naive search for a set of surface kernel generators in this group yields elements
$g_1$ and $g_2$ with orders 2 and 4 such that $(g_1g_2)^{-1}$ has
order 16.  There are four relations among these generators: 
\begin{displaymath}
g_1^2, \, g_2^4,  \, (g_2^{-1} g_1)^2 g_2^2 g_1 g_2 g_1 g_2^{-1}, \, (g_2 g_1)^2 g_2^{-1} (g_1 g_2)^2 g_1 g_2^{-1} (g_1 g_2)^2 g_1 
\end{displaymath}

\textit{Step 1.}
We use \texttt{Magma} to compute the conjugacy classes and character
table of $G$.  There are 16 conjugacy classes.  For convenience, write
$g_3 = g_1^{-1} g_2^{-1} g_1g_2$.  Then a list of representatives
of the conjugacy classes is

\[
\begin{array}{c}
\operatorname{Id},  \, g_3^4, \, g_2^2,\, g_1,\, g_3^2,\, g_2^2
  g_3^2,\, g_2^3,g_2,\\
g_1 g_2^2, \, g_3 g_2^2 g_3^4, \, g_3 g_2^2, \, g_3,g_2 g_1, \, g_2
  g_1 g_2^2 g_3^2, \, g_2 g_1 g_3^2, \, g_2 g_1 g_2^2
\end{array}
\]

Next we compute the character table.  The irreducible characters are given below
by their values on the sixteen conjugacy classes.
\begin{align*}
\chi_{1} & = ( 1 , 1 , 1 , 1 , 1 , 1 , 1 , 1 , 1 , 1 , 1 , 1 , 1 , 1 , 1 , 1)\\
\chi_{2} & = (1 , 1 , 1 , -1 , 1 , 1 , -1 , -1 , -1 , 1 , 1 , 1 , 1 , 1 , 1 , 1)\\
\chi_{3} & = (1 , 1 , 1 , 1 , 1 , 1 , -1 , -1 , 1 , 1 , 1 , 1 , -1 , -1 , -1 , -1)\\
\chi_{4} & = (1 , 1 , 1 , -1 , 1 , 1 , 1 , 1 , -1 , 1 , 1 , 1 , -1 , -1 , -1 , -1)\\
\chi_{5} & = (1 , 1 , -1 , 1 , 1 , -1 , -i , i , -1 , -1 , -1 , 1 , i , -i , i , -i)\\
\chi_{6} & = (1 , 1 , -1 , 1 , 1 , -1 , i , -i , -1 , -1 , -1 , 1 , -i , i , -i , i)\\
\chi_{7} & = (1 , 1 , -1 , -1 , 1 , -1 , i , -i , 1 , -1 , -1 , 1 , i , -i , i , -i)\\
\chi_{8} & = (1 , 1 , -1 , -1 , 1 , -1 , -i , i , 1 , -1 , -1 , 1 , -i , i , -i , i)\\
\end{align*}
\begin{align*}
\chi_{9} & = (2 , 2 , -2 , 0 , 2 , -2 , 0 , 0 , 0 , 2 , 2 , -2 , 0 , 0 , 0 , 0)\\
\chi_{10} & = (2 , 2 , 2 , 0 , 2 , 2 , 0 , 0 , 0 , -2 , -2 , -2 , 0 , 0 , 0 , 0)\\
\chi_{11} & = (2 , 2 , 2 , 0 , -2 , -2 , 0 , 0 , 0 , 0 , 0 , 0 ,-\sqrt{2} , \sqrt{2}, \sqrt{2}, -\sqrt{2})\\
\chi_{12} & = (2 , 2 , -2 , 0 , -2 , 2 , 0 , 0 , 0 , 0 , 0 , 0 , \sqrt{2}i , \sqrt{2}i , -\sqrt{2}i , -\sqrt{2}i)\\
\chi_{13} & = (2 , 2 , 2 , 0 , -2 , -2 , 0 , 0 , 0 , 0 , 0 , 0 , \sqrt{2}, -\sqrt{2}, -\sqrt{2}, \sqrt{2})\\
\chi_{14} & = (2 , 2 , -2 , 0 , -2 , 2 , 0 , 0 , 0 , 0 , 0 , 0 , -\sqrt{2}i , -\sqrt{2}i, \sqrt{2}i , \sqrt{2}i)\\
\chi_{15} & = (4 , -4 , 0 , 0 , 0 , 0 , 0 , 0 , 0 , -\sqrt{8}i , \sqrt{8}i  , 0 , 0 , 0 , 0 , 0)\\
\chi_{16} & = (4 , -4 , 0 , 0 , 0 , 0 , 0 , 0 , 0 , \sqrt{8}i  , -\sqrt{8}i  , 0 , 0 , 0 , 0 , 0)
\end{align*}

\textit{Step 2.}   Let $V_i$ be the irreducible $G$-module with character $\chi_i$ given
by the table above.  For any $G$-module $V$, let $V\cong
\bigoplus_{i=1}^{r} V_{i}^{\oplus m_i}$ be its decomposition into
irreducible $G$-modules.  

We use the Eichler trace formula in
\texttt{Magma} to compute these multiplicities $m_i$ for
several relevant $G$-modules.  Let $S=\mathbb{C}[x_0,x_1,x_2,x_3,x_4,x_5,x_6]$, and let $S_d$ denote polynomials of degree $d$.  Let
$I_d$ be the kernel defined by
\begin{displaymath}
0 \rightarrow I_d \rightarrow S_d \rightarrow \Gamma(\omega_X^{\otimes d}) \rightarrow 0.
\end{displaymath}
Then we have 
\begin{align*}
S_1 \cong \Gamma(\omega_X) & \cong V_8 \oplus V_{14} \oplus V_{15} \\
I_2 & \cong V_{3} \oplus V_{5} \oplus V_{10} \oplus V_{11} \oplus V_{16} \\
S_2 & \cong V_{3}^{\oplus 2} \oplus V_{5} \oplus V_{6} \oplus V_{10}^{\oplus 2} \oplus V_{11}^{\oplus 2} \oplus V_{13} \oplus V_{14} \oplus V_{15} \oplus V_{16}^{\oplus 2} \\
\Gamma(\omega_X^{\otimes 2})& \cong V_{3} \oplus V_{6} \oplus V_{10} \oplus V_{11} \oplus V_{13} \oplus V_{14} \oplus V_{15} \oplus V_{16}
\end{align*}

We use \texttt{GAP} to obtain matrix representatives of a $G$ action
with character equal to the character of the $G$ action on $S_1$.
Such a representation is obtained by mapping the generators $g_1$ and
$g_2$ to the matrices below.  
\begin{displaymath}
\left[ \begin{array}{rrrrrrr}
-1 & 0 & 0 & 0 & 0 & 0 & 0\\
0 & 0 & 1 & 0 & 0 & 0 & 0\\
0 & 1 & 0 & 0 & 0 & 0 & 0\\
0 & 0 & 0 & 0 & 1 & 0 & 0\\
0 & 0 & 0 & 1 & 0 & 0 & 0\\
0 & 0 & 0 & 0 & 0 & 0 & 1\\
0 & 0 & 0 & 0 & 0 & 1 & 0
\end{array} \right], \qquad
\left[ \begin{array}{rrrrrrr}
\zeta_{8}^2& 0 & 0 & 0 & 0 & 0 & 0\\
0 & 0 & -\zeta_{8} & 0 & 0 & 0 & 0\\
0 & -\zeta_{8}^3 & 0 & 0 & 0 & 0 & 0\\
0 & 0 & 0 & 0 & 0 & 0 & 1\\
0 & 0 & 0 & 0 & 0 & -\zeta_{8} & 0\\
0 & 0 & 0 & 0 &-\zeta_{8}^3 & 0 & 0\\
0 & 0 & 0 & 1 & 0 & 0 & 0
\end{array} \right]
\end{displaymath}
The decomposition of $S_1$ as a sum of three irreducible
$G$-modules gives rise to the block diagonal form of these matrices.  

\textit{Step 4.}  We use the projection
formula in \texttt{Magma} to decompose the $G$-module of quadrics $S_2$ into its
isotypical components.  When an isotypical component has multiplicity
greater than 1, we (noncanonically) choose ordered bases so that the $G$ action is
given by the same matrices on each ordered basis.  
\[
\begin{array}{lllll}
S_{2,3} & \cong & V_3^{\oplus 2} & = & \langle x_0^2\rangle
                                       \oplus \langle x_1x_2 \rangle \\
S_{2,5} & \cong & V_{5} & = & \langle x_3x_4 -\zeta_{8}x_5x_6 \rangle \\
S_{2,10} & \cong & V_{10}^{\oplus 2} &= & \langle x_1^2,x_2^2\rangle \oplus \langle x_3x_6+ix_4x_5, ix_3x_6+x_4x_5\rangle \\
S_{2,11} & \cong & V_{11}^{\oplus 2} & = & \langle x_0x_1, x_0x_2\rangle \oplus 
\langle x_3^2 +\zeta_{8}^3 x_5^2, -x_4^2 -\zeta_{8}^3x_6^2\rangle \\
S_{2,16} & \cong & V_{16}^{\oplus 2} & = &\langle x_0x_3,x_0x_4, x_0x_5,x_0x_6\rangle \oplus 
\langle -\zeta_{8} x_2 x_6,\zeta_{8} x_1 x_5,-x_2 x_4,x_1 x_3\rangle
\end{array}
\]

The first isotypical subspace yields a polynomial of the form $c_1
x_0^2 + c_2 x_1x_2$.  We may assume that $c_1$ and $c_2$ are nonzero,
scale $x_0$ to make $c_1=c_2$, and then divide by $c_1$ to obtain the
polynomial $x_0^2+x_1 x_2$. 

The second isotypical subspace yields the polynomial $x_3 x_4-\zeta_{8} x_5 x_6$.

The third isotypical subspace yields polynomials of the form $c_3 x_0 x_1+c_4 (x_3^2+\zeta_{8}^3 x_5^2)$ and 
$c_3x_0 x_2+c_4 (-x_4^2-\zeta_{8}^3 x_6^2)$.  We assume that $c_3$ and $c_4$ are nonzero, scale $x_1,x_2$ to make
$c_3=c_4$, then divide by $c_3$ and $c_4$.  

In the remaining isotypical subspaces no further scaling is possible, and hence we are left with two
undetermined coefficients $c_6$ and $c_8$.

Thus, a Riemann surface in this family has an ideal of the form 
\begin{displaymath}
\begin{array}{l}
x_0^2+x_1 x_2 \\
x_3 x_4-\zeta_{8} x_5 x_6\\
x_1^2+x_3 x_6+i x_4 x_5\\
x_2^2+i x_3 x_6+x_4 x_5 \\
x_0 x_1+c_6 (x_3^2+\zeta_{8}^3 x_5^2)\\
x_0 x_2+c_6 (-x_4^2-\zeta_{8}^3 x_6^2)\\
x_0 x_3+c_8 (-\zeta_{8} x_2 x_6)\\
x_0 x_4+c_8 (\zeta_{8} x_1 x_5)\\
x_0 x_5+c_8 (-x_2 x_4)\\
x_0 x_6+c_8 (x_1 x_3)
\end{array}
\end{displaymath}

\textit{Step 5.}  To find values of the coefficients $c_6,c_6$ that yield a smooth curve, we partially compute a flattening stratification.  
Begin Buchberger's algorithm.  We compute the S-pair reductions between
the generators and find  that
\begin{align*}
S(f_1,f_6) & \rightarrow (c_6 c_8 + \zeta_8^{-1})x_1x_4x_5+ \cdots \\
S(f_1,f_9) & \rightarrow ( c_8^2 - \zeta_8^{-1} ) x_1 x_2x_5 + \cdots
\end{align*}

Therefore, in Buchberger's algorithm, these polynomials will be added to the
Gr\"obner basis.  This suggests that we study the locus given by the equations
$c_8^2 -\zeta_8^{-1}=0$ and $c_6 c_8 + \zeta_8^{-1} = 0$ as an interesting stratum in the
flattening stratification.

We check in \texttt{Magma} that the values $c_6 = \zeta_{16}^{7}$ and $c_8 = \zeta_{16}^{-1}$ yield a smooth
genus 7 curve in $\Pro^6$ with the desired automorphism group.

From these equations, we can compute the Betti table of this ideal:
\begin{center}
\begin{tabular}{|c|c|c|c|c|c|}
\hline 1 & & & & & \\ \hline
& 10 & 16 & 3 & &  \\ \hline
&& 3 & 16 & 10 & \\ \hline
&&&&& 1 \\ \hline
\end{tabular}
\end{center}

Schreyer has classified Betti tables of genus 7 canonical curves in \cite{Schreyer1986}.
This Betti table implies that the curve is tetragonal (there exists a
degree 4 morphism $C \rightarrow
\Pro^1$) but not trigonal or hyperelliptic, and it has no degree
6 morphism $C \rightarrow \Pro^2$.

\section{Results} \label{results section}

This project had two goals.  The first goal was to establish that the
heuristics described in  Section \ref{heuristics section} allow us to run a variant of the main 
algorithm to completion for genus $4 \leq g \leq 7$ Riemann surfaces with large 
automorphism groups.  To this end, for each 
Riemann surface from Table 4 of \cite{MSSV}, the website
\cite{mywebpage} contains a link  to a calculation 
where a variant of the main algorithm is used to produce equations.  

The surface kernel generators needed to begin the algorithm were generally obtained by a naive search through the triples
or quadruples in the groups listed in Table 4 of \cite{MSSV}.  However, my \texttt{Magma}
code also includes functions allowing the user to input
surface kernel generators from any type of group, or to put in matrix
surface kernel generators with the desired representation on
$\Gamma(\omega_X)$. We note that Breuer's data has been recently extended and
republished by Paulhus \cite{Paulhus}, and Conder's data is available
online \cite{Conder}, so these sources could be used instead.

The equations obtained depend strongly on the matrix generators of the representation $\Aut(X)$ on 
$\Gamma(\omega_X)$.  I generally obtained these matrices from
\texttt{Magma}, \texttt{GAP}, the papers \cites{KKg34,KKg5}, or
\cite{Breuer}*{Appendix B}, and thus had little control over this step.  Indeed, in a few cases, the
resulting equations are almost comically bad; for an example
of this, compare my equations at \cite{mywebpage} for the genus 7 curve with 504 automorphisms 
to Macbeath's equations for this curve.  Given this, it is perhaps surprising that in most cases, the algorithm produces
reasonable equations (i.e., polynomials supported on a small number of
monomials 
with small coefficients).  

The second goal of this project was to create a reference that would 
contain the most useful information about the equations and automorphisms 
of these curves.  Thus, in this section, I print the best equations and automorphisms 
that I know, whether these were found in the literature or by the main algorithm.  
Many of the equations for genus $4 \leq g \leq 6$ are classical, and references are given whenever possible.  
However, the matrix surface kernel generators are not always equally easy to find.  The equations for the genus 
7 curves are almost all new, as are most of the 1-parameter families on the website \cite{mywebpage}.

\subsection{Description of the tables}
In the following tables I give equations for the Riemann surfaces 
of genus $4 \leq g \leq 7$ with large automorphism groups that 
are unique in moduli ($\delta=0$ in the notation of   Table 4 of \cite{MSSV}). 
The 1-parameter families ($\delta=1$) are not printed here but can be found
on the website \cite{mywebpage}.  I order the examples the same way they appear in \cite{MSSV}.  

For hyperelliptic Riemann surfaces, I give an equation of the form $y^2 = f(x)$.  Many of these 
are classically known, and all of them can be found in \cite{Shaska}.  

For plane quintics in genus 6, we give the plane quintic and surface kernel generators in $\GL(3,\mathbb{C})$.  
The canonical ideal and $G$ action can be easily computed from this data.  

For nonhyperelliptic curves that are not plane quintics, we print equations of the canonical ideals and surface kernel generators as
elements of $\GL(g,\mathbb{C})$.  Whenever such a matrix $M \in \GL(g,\mathbb{C})$ is sufficiently sparse, 
I frequently write the product $ M [x_0,\ldots, x_{g-1}]^{t}$ to save space.   

For the cyclic trigonal equations, I 
also print a cyclic trigonal equation, that is, one of the form 
$y^3 = \prod_{i=1}^{d_1} (x-\alpha_i) \prod_{i=1}^{d_2} (x-\beta_i)$, following the 
notation of \cite{AchterPries}*{Section 2.5} (where cyclic trigonal curves are also called \emph{trielliptic}).

Throughout the tables below, canonical ideals are shown in the polynomial ring $\mathbb{C}[x_0,\ldots,x_{g-1}]$.  
The symbol $\zeta_n$ denotes $e^{2 \pi i/n}$, and we write $i$ for $\zeta_4$.

\subsection{Genus 4}
In genus 4, every Riemann surface is either hyperelliptic or trigonal.  
Of the nine entries in 
Table 4 of \cite{MSSV}, four are 
hyperelliptic, four are cyclic trigonal, and one is general trigonal.

Note: the Riemann surface with automorphism group $(120,34) = S_5$ is known 
as Bring's curve.  Its best-known embedding is in $\Pro^5$, with equations $\sum_{i=0}^{4} x_i$, 
$\sum_{i=0}^{4} x_i^2$, $\sum_{i=0}^{4} x_i^3$.

\noindent Genus 4, Locus 1:  Group (120,34) = $S_5$, signature (2,4,5), general trigonal\\
\begin{tabular}{ll} Ideal:  & $x_0^2 + x_0 x_1 + x_1^2 - x_1x_2+ x_2^2 - x_2 x_3 + x_3^2$, \\ 
&  $x_0^2 x_1 + x_0 x_1^2 + x_1^2 x_2 - x_1x_2^2+ x_2^2 x_3 - x_2x_3^2$\end{tabular}\\
\begin{tabular}{ll} Maps: &$(x_0,x_1,x_2,x_3) \mapsto (-x_0,-x_1,-x_2,-x_2+x_3)$,\\
&  $(x_0,x_1,x_2,x_3) \mapsto (x_0+x_1,-x_0-x_2,-x_0-x_3,-x_3)$\end{tabular}\\
\mbox{} \\
Genus 4, Locus 2: Group (72,42), signature (2,3,12), cyclic trigonal\\
\begin{tabular}{l} Trigonal equation:  $y^3 = x(x^4-1)$\end{tabular}\\
\begin{tabular}{ll} Ideal:  & $x_1x_3-x_2^2,$ \\ 
& $x_0^3 - x_1^2 x_2 +x_2 x_3^2  $\end{tabular} \\
\begin{tabular}{l} Maps: $\left[ \begin{array}{rrrr}  
-1 & 0 & 0 & 0 \\
0 & 0 & 0 & -i \\
0 & 0 & -1 & 0 \\
0 & i & 0 & 0
\end{array} \right]$, \qquad 
$\left[ \begin{array}{rrrr}  
-\zeta_6 & 0 & 0 & 0 \\
0 & -\frac{1}{2}\zeta_{12} & \frac{1}{2}\zeta_{3}& \frac{1}{2}\zeta_{12} \\
0 & \zeta_{12} & 0 & \zeta_{12} \\
0 & -\frac{1}{2}\zeta_{12}& -\frac{1}{2}\zeta_{3} & \frac{1}{2}\zeta_{12}
\end{array} \right]$, 
\end{tabular}\\
\mbox{} \\
Genus 4, Locus 3: Group (72,40), signature (2,4,6), cyclic trigonal\\
\begin{tabular}{l} Trigonal equation:  $y^3 = (x^3-1)^2(x^3+1)$\end{tabular}\\
\begin{tabular}{ll} Ideal:  &$x_0 x_3-x_1 x_2,$ \\ 
&  $x_1^3-x_0^3-x_3^3-x_2^3  $\end{tabular}\\
\begin{tabular}{ll} Maps: & $(x_0,x_1,x_2,x_3) \mapsto(-x_0,x_2,x_1,-x_3)$\\
& $(x_0,x_1,x_2,x_3) \mapsto (-x_2,\zeta_6^2x_0,\zeta_6x_3,x_1)$\end{tabular}\\
\mbox{} \\
Genus 4, Locus 4: Group (40,8), signature (2,4,10), hyperelliptic\\
\begin{tabular}{l} $y^2=x^{10}-1$ \end{tabular} \\
Genus 4, Locus 5: Group (36,12), signature (2,6,6), cyclic trigonal\\
\begin{tabular}{l}  Trigonal equation:  $y^3 = (x^3-1)(x^3+1)$\end{tabular}\\
\begin{tabular}{ll}  Ideal: & $x_1 x_3-x_2^2,$ \\ 
&  $x_0^3-x_3^3+x_1^3   $\end{tabular}\\
\begin{tabular}{ll}  Maps: &$(x_0,x_1,x_2,x_3) \mapsto(-x_0,\zeta_3 x_3,-x_2,-\zeta_6 x_1)$\\
& $(x_0,x_1,x_2,x_3) \mapsto (\zeta_3 x_0,-\zeta_3 x_3, \zeta_6 x_2,
  -x_1)$\end{tabular}\\
\newpage
\noindent Genus 4, Locus 6: Group (32,19), signature (2,4,16), hyperelliptic\\
\begin{tabular}{l} $y^2=x^9-x$ \end{tabular} \\
\mbox{} \\
Genus 4, Locus 7: Group (24,3), signature (3,4,6), hyperelliptic\\
\begin{tabular}{l} $y^2=x(x^4-1)(x^4+2i\sqrt{3}+1)$ \end{tabular} \\
\mbox{} \\
Genus 4, Locus 8: Group (18,2), signature (2,9,18), hyperelliptic\\
\begin{tabular}{l} $y^2=x^9-1$ \end{tabular} \\
\mbox{} \\
Genus 4, Locus 9:  Group (15,1), signature (3,5,15), cyclic trigonal\\
\begin{tabular}{l} Trigonal equation:  $y^3 = x^5-1$ \end{tabular}\\
\begin{tabular}{ll}  Ideal: & $x_1 x_3 -x_2^2,$ \\ 
& $x_0^3 - x_1^2 x_2 +x_3^3   $\end{tabular}\\
\begin{tabular}{ll}  Maps: & $(x_0,x_1,x_2,x_3) \mapsto (\zeta_{3}^{2} x_0,\zeta_{3} x_1, \zeta_{3}x_2, \zeta_{3} x_3)$\\
& $(x_0,x_1,x_2,x_3) \mapsto (\zeta_{5} x_0,\zeta_{5}^3 x_1,\zeta_{5}^2 x_2,\zeta_{5} x_3)$
\end{tabular}

\subsection{Genus 5}
Of the ten entries in Table 4 of \cite{MSSV}, five 
are hyperelliptic, and one is cyclic trigonal.  The remaining four are general, 
hence their canonical models are complete intersections of three quadrics.

\noindent Genus 5, Locus 1: Group (192,181), signature (2,3,8) \\
\begin{tabular}{ll} Ideal: &  Wiman, \cite{Wiman}:\\
&  $x_0^2+x_3^2+x_4^2$, \\ 
&  $x_1^2+x_3^2-x_4^2$\\
& $x_2^2+x_3 x_4$\end{tabular}\\
\begin{tabular}{ll} Maps: & $\left[\begin{array}{rrrrr}
0 & 0 & \frac{1}{2}(i+1) & 0 & 0 \\
  0 & -1 & 0 & 0 & 0\\
  1-i & 0 & 0 & 0 & 0\\
  0 & 0 & 0 & -\frac{1}{\sqrt{2}} & -\frac{i}{\sqrt{2}}\\
  0 & 0 & 0 & \frac{i}{\sqrt{2}} & \frac{1}{\sqrt{2}} 
  \end{array} \right]$,\\ & $\quad \left[\begin{array}{rrrrr}
0 & \zeta_8^{-1} & 0 & 0 & 0 \\
  0 & 0 & -\frac{1}{\sqrt{2}}  & 0 & 0\\
  -1-i & 0 & 0 & 0 & 0\\
  0 & 0 & 0 & \frac{1}{2}(i-1) & -\frac{1}{2}(i+1)\\
  0 & 0 & 0 & -\frac{1}{2}(i-1) & -\frac{1}{2}(i+1)
\end{array} \right]$ \end{tabular}\\
\mbox{} \\
Genus 5, Locus 2:  Group (160,234), signature (2,4,5) \\
\begin{tabular}{ll} Ideal: &  Wiman, \cite{Wiman}:\\
& $x_0^2+x_1^2+x_2^2+x_3^2+x_4^2$, \\ 
&  $x_0^2+\zeta_5 x_1^2+\zeta_5^2 x_2^2+\zeta_5^3 x_3^2+\zeta_5^4 x_4^2,$\\
& $\zeta_5^4 x_0^2+\zeta_5^3 x_1^2+\zeta_5^2 x_2^2+\zeta_5 x_3^2+x_4^2$\end{tabular}\\
\begin{tabular}{ll}  Maps: & $(x_0,x_1,x_2,x_3,x_4) \mapsto (-x_3,x_2,x_1,-x_0,-x_4)$,\\
&  $(x_0,x_1,x_2,x_3,x_4) \mapsto (-x_0,x_4,-x_3,x_2,-x_1)$\end{tabular}\\
\mbox{} \\
Genus 5, Locus 3: Group (120,35), signature (2,3,10), hyperelliptic\\
\begin{tabular}{l} $y^2=x^{11}+11x^6-x$\end{tabular}\\
\mbox{}\\
Genus 5, Locus 4:  Group (96,195), signature (2,4,6) \\
\begin{tabular}{ll} Ideal: &  Wiman, \cite{Wiman}:\\
& $x_0^2 + x_3^2 + x_4^2$, \\ 
& $x_1^2+\zeta_3 x_3^2+\zeta_3^2 x_4^2,$\\
& $x_2^2+\zeta_3^2 x_3^2+\zeta_3 x_4^2$\end{tabular}\\
\begin{tabular}{ll} Maps: & $(x_0,x_1,x_2,x_3,x_4) \mapsto (-x_2,-x_1,-x_0,\zeta_3^2 x_4,\zeta_3 x_3)$,\\
&  $(x_0,x_1,x_2,x_3,x_4) \mapsto
  (-x_0,x_2,-x_1,-x_4,x_3)$\end{tabular}\\
\newpage
\noindent Genus 5, Locus 5:  Group (64,32), signature (2,4,8) \\
\begin{tabular}{ll} Ideal: &  Wiman, \cite{Wiman}:\\
& $x_0^2+x_1^2+x_2^2+x_3^2+x_4^2$, \\ 
& $x_0^2+i x_1^2-x_2^2-i x_3^2,$\\
& $x_0^2-x_1^2+x_2^2-x_3^2$\end{tabular} \\
\begin{tabular}{ll} Maps: &$(x_0,x_1,x_2,x_3,x_4) \mapsto (-x_0,x_1,-x_2,-x_3,-x_4)$,\\
&  $(x_0,x_1,x_2,x_3,x_4) \mapsto (ix_1,-ix_2,ix_3,-ix_0,ix_4)$\end{tabular}\\
\mbox{} \\
Genus 5, Locus 6: Group (48,14), signature (2,4,12), hyperelliptic\\
\begin{tabular}{l} $y^2=x^{12}-1$\end{tabular}\\
\mbox{}\\
Genus 5, Locus 7: Group (48,30), signature (3,4,4), hyperelliptic\\
\begin{tabular}{l} $y^2=x^{12}-33x^8-33x^{4}+1$\end{tabular}\\
\mbox{}\\
Genus 5, Locus 8: Group (40,5), signature (2,4,20), hyperelliptic\\
\begin{tabular}{l} $y^2=x^{11}-x$\end{tabular}\\
\mbox{}\\
Genus 5, Locus 9: Group (30,2), signature (2,6,15), cyclic trigonal \\
\begin{tabular}{l}  Trigonal equation: $y^3=(x^5-1)x^2$ \end{tabular}\\
\begin{tabular}{ll} Ideal: & $x_0 x_3-x_1 x_2$, \,  $x_0 x_4-x_1 x_3$,
                             \, $x_2 x_4-x_3^2,$\\
& $x_0^2 x_1-x_3 x_4^2+x_2^3,$\\
& $x_0 x_1^2-x_4^3+x_2^2 x_3$ \end{tabular}\\
\begin{tabular}{ll} Maps: & $(x_0,x_1,x_2,x_3,x_4) \mapsto (\zeta_5 x_1,\zeta_5^4 x_0,-\zeta_5^2 x_4,-x_3,-\zeta_5^3 x_2)$,\\
&  $(x_0,x_1,x_2,x_3,x_4) \mapsto (\zeta_{15}^{14} x_1,\zeta_{15}^{11} x_0,-\zeta_{15}^{13} x_4,-\zeta_{15}^{10} x_3,-\zeta_{15}^7x_2)$\end{tabular}\\
\mbox{} \\
Genus 5, Locus 10: Group (22,2), signature (2,11,22), hyperelliptic\\
\begin{tabular}{l} $y^2=x^{11}-1$\end{tabular}

\subsection{Genus 6}
Table 4 in \cite{MSSV} contains eleven entries for genus 6 Riemann surfaces 
with large automorphism groups and no moduli ($\delta=0$). Of these, four 
are hyperelliptic, three are cyclic trigonal, and three are plane quintics; only one is general.  

For the plane quintics, we give the plane quintic equation in the variables $y_0,y_1,y_2$, and surface kernel generators
acting on the plane.  The canonical model of a plane quintic lies on the Veronese surface, and 
 the multiples of the quintic by $y_0,y_1,y_2$ may be encoded as cubics in $x_0,\ldots,x_5$.  \\

\noindent Genus 6, Locus 1: Group (150,5), signature (2,3,10), plane quintic \\
\begin{tabular}{l} Plane quintic equation:  $y_0^5+y_1^5+y_2^5$\end{tabular}\\ 
\begin{tabular}{ll} Maps: & $(y_0,y_1,y_2) \mapsto (-\zeta_5^3 y_1, -\zeta_5^2 y_0, -y_2)$,\\
& $(y_0,y_1,y_2) \mapsto (-\zeta_5^3 y_1, -\zeta_5^2 y_0, -y_2)$
\end{tabular}\\
\mbox{} \\
\noindent Genus 6, Locus 2: Group $(120,34) = S_5$, signature (2,4,6) \\
\begin{tabular}{ll} Ideal: & Inoue and Kato,\cite{InoueKato}:\\
 &$-x_0x_2 + x_1x_2 - x_0x_3 + x_1x_4,$\\
& $ -x_0x_1 + x_1x_2 - x_0x_3 + x_2x_5,$ \\
& $ -x_0x_1 - x_0x_2 - 2x_0x_3 - x_3x_4 - x_3x_5,$ \\
& $ -x_0x_1 - x_0x_2 - x_0x_3 - x_1x_4 - x_3x_4 - x_4x_5,$ \\
& $ -x_0x_1 - x_0x_2 - x_0x_3 - x_2x_5 - x_3x_5 - x_4x_5$, \\
& $2 (\sum_{i=1}^{6} x_i^2) + x_0 x_1 + x_0 x_2 + x_1 x_2 + 2 x_1 x_3 + 2 x_2 x_3 $\\
& $\mbox{        }+2 x_0 x_4 + 2 x_2 x_4 + x_3 x_4 + 2 x_0 x_5 + 2 x_1 x_5 + x_3 x_5 + x_2 x_5 $ \end{tabular}\\
\begin{tabular}{l} Maps: $\left[\begin{array}{rrrrrr}
0 & 0 & 0 & 0 & -1 & 0 \\
0 & -1 & 1 & 0 & -1 & 1 \\
-1 & 0 & 0 & 0 & 0 & 1 \\
-1 & 0 & 1 & -1 & 0 & 1 \\
-1 & 0 & 0 & 0 & 0 & 0 \\
0 & 0 & 1 & 0 & -1 & 0 
\end{array} \right]$, \quad $\left[\begin{array}{rrrrrr}
0 & 0 & 1 & 0 & -1 & 0 \\
-1 & 0 & 1 & -1 & 0 & 1\\
-1 & 0 & 0 & 0 & 0 & 0 \\
-1 & 0 & 0 & 0 & 0 & 1 \\
0 & 0 & 0 & 0 & -1 & 0 \\
0 & -1 & 1 & 0 & -1 & 0 
\end{array} \right]$
\end{tabular}\\
\mbox{} \\
Genus 6, Locus 3:  Group (72,15), signature (2,4,9), cyclic trigonal
\\
\begin{tabular}{l} Trigonal equation:  $y^3 = (x^4-2 \sqrt{3} i x^2+1)(x^4+2 \sqrt{3}ix^2 + 1)^2$\end{tabular}\\
\begin{tabular}{ll}Ideal: & $x_0 x_2 - x_1^2$, $x_0 x_4 - x_1 x_3 $, $x_0 x_5-x_1 x_4 $,\\
&  $ x_1 x_4 - x_2 x_3$, $x_1x_5-x_2x_4 $, $x_3x_5-x_4^2 $,\\
& $x_0^3+(4 \zeta_6-2) x_0^2 x_2+x_0 x_2^2+x_3^3+(-4 \zeta_6+2) x_3^2 x_5+x_3 x_5^2$\\
& $x_0^2 x_1+(4 \zeta_6-2) x_0 x_1 x_2+x_1 x_2^2+x_3^2 x_4+(-4 \zeta_6+2) x_3 x_4 x_5+x_4 x_5^2,$\\
& $x_0^2 x_2+(4 \zeta_6-2) x_0 x_2^2+x_2^3+x_3^2 x_5+(-4 \zeta_6+2)
  x_3 x_5^2+x_5^3$\end{tabular} \\
\begin{tabular}{ll} Maps: & $\left[\begin{array}{rrrrrr}
0 & 0 & 0 & -\frac{1}{2}\zeta_9 & \frac{1}{2}\zeta_{36}^{13} & \frac{1}{2}\zeta_9\\
0 & 0 & 0 & -\zeta_{36}^{13} & 0 & -\zeta_{36}^{13}\\
0 & 0 & 0 & \frac{1}{2}\zeta_9 & \frac{1}{2}\zeta_{36}^{13} & -\frac{1}{2}\zeta_9\\
\frac{1}{2}\zeta_{36}^{14} & \frac{1}{2}\zeta_{36}^5 & -\frac{1}{2}\zeta_{36}^{14} & 0 & 0 & 0\\
-\zeta_{36}^5 & 0 & -\zeta_{36}^5 & 0 & 0 & 0\\
-\frac{1}{2}\zeta_{36}^{14} & \frac{1}{2}\zeta_{36}^5 & \frac{1}{2}\zeta_{36}^{14} & 0 & 0 & 0
\end{array} \right],$\\
& $(x_0,\ldots,x_5) \mapsto (\zeta_{12}x_3,\zeta_3 x_4, -\zeta_{12} x_5, \zeta_{12}^5x_0,\zeta_3^2 x_1,-\zeta_{12}^5x_2) $\end{tabular}\\
\mbox{} \\
Genus 6, Locus 4:  Group (56,7), signature (2,4,14), hyperelliptic\\
\begin{tabular}{l} $y^2=x^{14}-1$ \end{tabular}\\
\mbox{}\\
Genus 6, Locus 5:  Group (48,6), signature (2,4,24), hyperelliptic\\
\begin{tabular}{l} $y^2=x^{13}-x$ \end{tabular}\\
\mbox{} \\
Genus 6, Locus 6:  Group (48,29), signature (2,6,8), hyperelliptic \\
\begin{tabular}{l} $y^2=x(x^4-1)(x^{8}+14x^4+1)$ \end{tabular}\\
\mbox{} \\
Genus 6, Locus 7:  Group (48,15), signature (2,6,8), cyclic trigonal \\
\begin{tabular}{l} Trigonal equation:  $y^3 = (x^4-1)^2(x^4+1)$\end{tabular}\\
\begin{tabular}{ll} Ideal: & $x_0 x_2 - x_1^2$, $x_0 x_4 - x_1 x_3 $, $x_0 x_5-x_1 x_4 $,\\
&  $ x_1 x_4 - x_2 x_3$, $x_1x_5-x_2x_4 $, $x_3x_5-x_4^2 $,\\
&  $x_0 x_1^2 - x_2^3 - x_3 x_4^2 - x_5^3, $\\
&  $x_0^2 x_1 - x_1 x_2^2 - x_3^2 x_4 - x_4 x_5^2,$\\
&  $x_0^3 - x_1^2 x_2 - x_3^3 - x_4^2 x_5$\end{tabular}\\
\begin{tabular}{ll} Maps: & $(x_0,\ldots,x_5) \mapsto (\zeta_8^3 x_5, -i x_4, \zeta_8 x_3, -\zeta_8^3 x_2, i x_1, -\zeta_8 x_0)$,\\
&  $(x_0,\ldots,x_5)\mapsto (-\zeta_6 x_2,\zeta_6 x_1, -\zeta_6 x_0, -\zeta_3 x_5,\zeta_3 x_4,-\zeta_3 x_3)$\end{tabular}\\
\mbox{} \\
Genus 6, Locus 8:  Group (39,1), signature (3,3,13), plane quintic \\
\begin{tabular}{l} Plane quintic equation:  $y_0^4y_1 + y_1^4y_2 + y_2^4y_0$\end{tabular}\\ 
\begin{tabular}{ll} Maps: & $(y_0,y_1,y_2) \mapsto (\zeta_{13}^4 y_1, \zeta_{13}^{10} y_2, \zeta_{13}^{12} y_0)$,\\
&  $(y_0,y_1,y_2) \mapsto (\zeta_{13}^8 y_2, \zeta_{13}^{7} y_0, \zeta_{13}^{11} y_1)$\end{tabular}\\
\mbox{} \\
Genus 6, Locus 9:  Group (30,1), signature (2,10,15), plane quintic \\
\begin{tabular}{l} Plane quintic equation:  $y_0^5+y_1^4 y_2+2 \zeta_5 y_1^3 y_2^2+2 \zeta_5^2 y_1^2 y_2^3+\zeta_5^3 y_1 y_2^4$\end{tabular}\\ 
\begin{tabular}{ll} Maps: & $(y_0,y_1,y_2) \mapsto (y_0,\zeta_5 y_2, \zeta_5^4 y_1)$,\\
&  $(y_0,y_1,y_2) \mapsto (\zeta_5^3 y_0, -\zeta_5y_1-\zeta_5^2y_2,\zeta_5)$\end{tabular}\\
\mbox{} \\
Genus 6, Locus 10:  Group (26,2), signature (2,13,26), hyperelliptic\\
\begin{tabular}{l} $y^2=x^{13}-1$ \end{tabular}\\
Genus 6, Locus 11:  Group (21,2), signature (3,7,21), cyclic trigonal \\
\begin{tabular}{l} Trigonal equation:  $y^3 = x^7-1$\end{tabular} \\
\begin{tabular}{ll} Ideal:  &$x_0 x_3-x_1 x_2 $, $x_0 x_4-x_1 x_3 $, $x_0 x_5-x_1 x_4 $, \\
& $x_2 x_4-x_3^2 $, $x_2 x_5-x_3 x_4 $, $x_3 x_5-x_4^2 $, \\
& $x_0^3-x_3 x_5^2 +x_2^3 $,\\
& $x_0^2 x_1-x_4 x_5^2 + x_2^2 x_3$,\\
& $x_0 x_1^2-x_5^3+x_2^2 x_4 $\end{tabular} \\
\begin{tabular}{ll} Maps:& $(x_0,\ldots,x_5) \mapsto (x_0,\zeta_7 x_1, x_2, \zeta_7 x_3,\zeta_7^2 x_4 \zeta_7^3 x_5)$,\\
&  $(x_0,\ldots,x_5)\mapsto (\zeta_3^2 x_0, \zeta_3^2 x_1, \zeta_3 x_2, \zeta_3 x_3, \zeta_3 x_4, \zeta_3 x_5)$
\end{tabular}

\subsection{Genus 7}
Of the thirteen entries in Table 4 of \cite{MSSV} for genus 7 curves, three 
are hyperelliptic and two are cyclic trigonal.  After computing the
canonical equations of the nonhyperelliptic Riemann surfaces, we can
compute the Betti tables of these ideals and use the results of \cite{Schreyer1986} to classify the curve
as having a $g^{1}_{4}$, $g^{2}_{6}$, $g^{1}_{3}$, or none of these.  

\noindent Genus 7, Locus 1: Group $(504,156)$, signature (2,3,7) \\
\begin{tabular}{ll} Ideal: & Macbeath,\cite{Macbeath}: \\
&  $ x_0^2+x_1^2+x_2^2+x_3^2+x_4^2+x_5^2+x_6^2,$ \\
&  $ x_0^2+\zeta_{7} x_1^2+\zeta_{7}^2 x_2^2+\zeta_{7}^3 x_3^2+\zeta_{7}^4 x_4^2+\zeta_{7}^5 x_5^2+\zeta_{7}^6 x_6^2,$ \\
&  $x_0^2+\zeta_{7}^{-1} x_1^2+\zeta_{7}^{-2} x_2^2+\zeta_{7}^{-3} x_3^2+\zeta_{7}^{-4} x_4^2+\zeta_{7}^{-5} x_5^2+\zeta_{7}^{-6} x_6^2,$ \\
&  $(\zeta_{7}^{-3}-\zeta_{7}^3) x_0 x_6-(\zeta_{7}^{-2}-\zeta_{7}^2) x_1 x_4+(\zeta_{7}-\zeta_{7}^{-1}) x_3 x_5,$ \\
&  $(\zeta_{7}^{-3}-\zeta_{7}^3) x_1 x_0-(\zeta_{7}^{-2}-\zeta_{7}^2) x_2 x_5+(\zeta_{7}-\zeta_{7}^{-1}) x_4 x_6,$ \\
&  $(\zeta_{7}^{-3}-\zeta_{7}^3) x_2 x_1-(\zeta_{7}^{-2}-\zeta_{7}^2) x_3 x_6+(\zeta_{7}-\zeta_{7}^{-1}) x_5 x_0,$ \\
&  $(\zeta_{7}^{-3}-\zeta_{7}^3) x_3 x_2-(\zeta_{7}^{-2}-\zeta_{7}^2) x_4 x_0+(\zeta_{7}-\zeta_{7}^{-1}) x_6 x_1,$ \\
&  $(\zeta_{7}^{-3}-\zeta_{7}^3) x_4 x_3-(\zeta_{7}^{-2}-\zeta_{7}^2) x_5 x_1+(\zeta_{7}-\zeta_{7}^{-1}) x_0 x_2,$ \\
&  $(\zeta_{7}^{-3}-\zeta_{7}^3) x_5 x_4-(\zeta_{7}^{-2}-\zeta_{7}^2) x_6 x_2+(\zeta_{7}-\zeta_{7}^{-1}) x_1 x_3,$ \\
&  $(\zeta_{7}^{-3}-\zeta_{7}^3) x_6 x_5-(\zeta_{7}^{-2}-\zeta_{7}^2) x_0 x_3+(\zeta_{7}-\zeta_{7}^{-1}) x_2 x_4$ \end{tabular}\\
\begin{tabular}{ll} Maps: & $(x_0,\ldots,x_6)\mapsto (x_0,-x_1,-x_2,-x_3,x_4,x_5,-x_6 )$,\\
& $\left[\begin{array}{rrrrrrr}
0 & \frac{1}{2} & \frac{1}{2} & -\frac{1}{2} & 0 & -\frac{1}{2} & 0 \\
-\frac{1}{2} & -\frac{1}{2} & \frac{1}{2} & 0 & -\frac{1}{2} & 0 & 0 \\
\frac{1}{2} & -\frac{1}{2} & 0 & -\frac{1}{2} & 0 & 0 & -\frac{1}{2}\\
-\frac{1}{2} & 0 & -\frac{1}{2} & 0 & 0 & -\frac{1}{2} & -\frac{1}{2} \\
0 & \frac{1}{2} & 0 & 0 & -\frac{1}{2} & \frac{1}{2} & -\frac{1}{2} \\
\frac{1}{2} & 0 & 0 & \frac{1}{2} & -\frac{1}{2} & -\frac{1}{2} & 0 \\
0 & 0 & \frac{1}{2} & \frac{1}{2} &\frac{1}{2} & 0 & -\frac{1}{2}
\end{array} \right]$
\end{tabular}\\
\mbox{} \\
Genus 7, Locus 2: Group $(144,127)$, signature (2,3,12). Has $g^{2}_{6}$ \\
\begin{tabular}{ll} Ideal: &$x_0^2 + x_3x_4 - \zeta_6x_3x_5 - \zeta_6x_5x_6,$ \\
&$2ix_1^2 + x_3x_4 + \zeta_6x_3x_5 + 2x_4x_6 - \zeta_6x_5x_6,$ \\
&$2ix_1x_2 + (-2\zeta_6 + 1)x_3x_4 + \zeta_6x_3x_5 - 2\zeta_6x_4x_6 + \zeta_6x_5x_6,$ \\
&$2ix_2^2 -x_3x_4 + (-\zeta_6 + 2)x_3x_5 + (2\zeta_6 - 2)x_4x_6 + (-\zeta_6 + 2)x_5x_6,$ \\
& $x_1x_3 - \zeta_6x_2x_6 + \zeta_{12}x_4^2 + (\zeta_{12}^3 - 2\zeta_{12})x_4x_5 + \zeta_{12}x_5^2,$ \\
& $x_1x_4 -\zeta_3x_2x_5 -x_3x_6 - x_6^2,$ \\
& $x_1x_5 - x_2x_5 + x_3^2 + (-\zeta_6 + 2)x_3x_6,$ \\
& $x_1x_6 + x_2x_6   -\zeta_{12}x_4^2 + \zeta_{12}x_4x_5 - \zeta_{12}x_5^2,$ \\
& $x_2x_3 - \zeta_6x_2x_6 + \zeta_{12}x_4^2,$ \\
& $x_2x_4 + (-\zeta_6 - 1)x_2x_5 + \zeta_6x_3^2 + 2\zeta_6x_3x_6 + \zeta_6 x_6^2$\end{tabular}\\

\begin{tabular}{ll}Maps: &$\left[\begin{array}{rrrrrrr}
-1 & 0 & 0 & 0 & 0 & 0 & 0 \\
0 & \zeta_{12}^{-1} & -\zeta_{12} &  0 & 0 & 0 & 0 \\
0 & -\zeta_{12} & -\zeta_{12}^{-1} & 0 & 0 & 0 & 0\\
0 & 0 & 0 & 0 & \zeta_6 & \zeta_6 & 0 \\
0 & 0 & 0 & -\zeta_3 & 0 & 0 & \zeta_3\\
0 & 0 & 0 & 0 & 0 & 0 & -\zeta_3 \\
0 & 0 & 0 & 0 & 0 & \zeta_6 & 0
\end{array} \right]$, \\
& $\left[\begin{array}{rrrrrrr}
\zeta_3 & 0 & 0 & 0 & 0 & 0 & 0 \\
0 & -1 & \zeta_3 & 0 & 0 & 0 & 0 \\
0 & \zeta_6 & 0 & 0 & 0 & 0 & 0 \\
0 & 0 & 0 & \zeta_3 & 0 & 0 & 1 \\
0 & 0 & 0 & 0 & \zeta_3 & 0 & 0 \\
0 & 0 & 0 & 0 & -\zeta_3 & -\zeta_6 & 0 \\
0 & 0 & 0 & 0 & 0 & 0 & 1
\end{array} \right]$\end{tabular}\\
\newpage
\noindent Genus 7, Locus 3: Group $(64,41)$, signature (2,3,16), tetragonal \\
\begin{tabular}{ll}Ideal: & $x_0^2+x_1 x_2,$ \\
& $x_3 x_4-\zeta_8 x_5 x_6,$ \\
& $x_1^2+x_3 x_6+i x_4 x_5,$ \\
& $x_2^2+i x_3 x_6+x_4 x_5,$ \\
& $x_0 x_1+\zeta_{16}^7 x_3^2-\zeta_{16}^5 x_5^2,$ \\
& $x_0 x_2-\zeta_{16}^7 x_4^2+\zeta_{16}^5 x_6^2,$ \\
& $x_0 x_3-\zeta_{16} x_2 x_6,$ \\
& $x_0 x_4+\zeta_{16} x_1 x_5,$ \\
& $x_0 x_5+\zeta_{16}^7 x_2 x_4,$ \\
& $x_0 x_6-\zeta_{16}^7 x_1 x_3$\end{tabular} \\
\begin{tabular}{ll}Maps:& $(x_0,\ldots,x_6) \mapsto  (-x_0,x_2,x_1,x_4,x_3,x_6,x_5)$, \\
& $(x_0,\ldots,x_6) \mapsto  (ix_0, -\zeta_{8}^3 x_2,-\zeta_{8}x_1,x_6,-\zeta_{8}^3x_5,-\zeta_{8} x_4,x_3)$\end{tabular}\\
\mbox{} \\
Genus 7, Locus 4: Group $(64,38)$, signature (2,4,16), hyperelliptic \\
\begin{tabular}{l} $y^2=x^{16}-1$ \end{tabular} \\
\mbox{}\\
Genus 7, Locus 5: Group $(56,4)$, signature (2,4,28), hyperelliptic \\
\begin{tabular}{l} $y^2=x^{15}-x$ \end{tabular} \\
\mbox{}\\
Genus 7, Locus 6: Group $(54,6)$, signature (2,6,9) \\
\begin{tabular}{ll}Ideal: & $x_1 x_6 + x_2 x_4 + x_3 x_5,$ \\
& $x_0^2-x_1 x_6 +\zeta_6 x_2 x_4 -\zeta_3 x_3 x_5,$ \\
& $x_1 x_4 + \zeta_3 x_2 x_5 - \zeta_6 x_3 x_6,$ \\
& $x_1 x_5 + \zeta_3 x_2 x_6 - \zeta_6 x_3 x_4,$ \\
& $x_0 x_1-\zeta_6 x_5^2-x_4 x_6,$ \\
& $x_0 x_2+x_6^2-\zeta_3 x_4 x_5,$ \\
& $x_0 x_3+\zeta_3 x_4^2+\zeta_6 x_5 x_6,$ \\
& $x_0 x_4-x_1^2+\zeta_3 x_2 x_3,$ \\
& $x_0 x_5-\zeta_3 x_2^2-\zeta_6 x_1 x_3,$ \\
& $x_0 x_6+\zeta_6 x_3^2+x_1 x_2$\end{tabular}\\
\begin{tabular}{ll}Maps: & $(x_0,\ldots,x_6) \mapsto  (-x_0,\zeta_9^5 x_6,\zeta_9^8 x_4,\zeta_9^2 x_5,\zeta_9 x_2,\zeta_9^7 x_3,\zeta_9^4 x_1 )$, \\
& $(x_0,\ldots,x_6) \mapsto  (\zeta_6 x_0, \zeta_3^2 x_4, \zeta_3^2
  x_5, \zeta_3^2 x_6, \zeta_3 x_3, \zeta_3 x_1, \zeta_3 x_2)$\end{tabular}\\
\mbox{} \\
Genus 7, Locus 7: Group $(54,6)$, signature (2,6,9) \\
\begin{tabular}{l}Complex conjugate of the previous curve \end{tabular}\\
\mbox{} \\

Genus 7, Locus 8:  Group $(54,3)$, signature (2,6,9) cyclic trigonal\\
\begin{tabular}{l} Trigonal equation: $y^3 = x^9-1$\end{tabular}\\
\begin{tabular}{ll} Ideal: &  $2\times 2$ minors of 
$\left[ \begin{array}{rrrrr} x_0 & x_2 & x_3 & x_4 & x_5 \\ x_1 & x_3 & x_4 & x_5 & x_6 \end{array} \right]$, and \\
& $ x_0^3-x_6^2+x_2^3$, \\
& $x_0^2 x_1-x_6^2 x_4+x_2^2 x_3,$ \\
& $x_0 x_1^2-x_6^2 x_5+x_2^2 x_4,$ \\
& $x_1^3-x_6^3+x_2^2 x_5$ \end{tabular}\\
\begin{tabular}{ll} Maps: & $(x_0,\ldots,x_6) \mapsto  (x_1,x_0,-x_6,-x_5,-x_4,-x_3,-x_2)$, \\
& $(x_0,\ldots,x_6) \mapsto  (\zeta_9 x_1, \zeta_9^2 x_0,-\zeta_9 x_6, -\zeta_9^2 x_5, -\zeta_9^3 x_4, -\zeta_9^4 x_3, -\zeta_9^5 x_2 )$\end{tabular}\\
\setlength{\textheight}{9.25in}
\newpage
\noindent Genus 7, Locus 9:  Group $(48,32)$, signature (3,4,6). Has $g^{2}_{6}$ \\
\begin{tabular}{l}Ideal: \\
$x_0^2+x_3 x_5 + \zeta_6 x_3 x_6 -\zeta_3 x_4 x_6,$ \\
$\sqrt{3}i (x_1 x_3 - x_2 x_4+ x_1 x_5) - x_1 x_6 + x_2 x_5 - x_2 x_6,$ \\
$\sqrt{3}i (2 x_1 x_4 - x_2 x_5 + x_2 x_6)+ x_1 x_5 - x_1 x_6,$ \\
$\sqrt{3}i (2 x_2 x_3 + x_1 x_6 - x_2 x_5)  + 3 x_1 x_5 + x_2 x_6,$ \\
$-3 (x_1^2+x_3 x_5 - x_3 x_6) + \sqrt{3}i (x_4 x_5 - x_4 x_6)+2 x_5^2 + 2 (\zeta_6 - 1) x_5 x_6 -2 \zeta_6 x_6^2,$ \\
$-3 (x_1 x_2 - x_4 x_5)+\sqrt{3}i (x_3 x_5 - x_3 x_6 + x_4 x_6-x_5^2) + 2 \zeta_6 x_5 x_6 - x_6^2,$ \\
$-3 (x_2^2-x_3 x_5 - x_4 x_6)+ \sqrt{3}i (x_3 x_6 -3 x_4 x_5) +2 (\zeta_6 + 1) x_5 x_6 + 2 (\zeta_6 - 1) x_6^2,$ \\
$-3 (2 x_4^2-x_3 x_5 + x_3 x_6) +\sqrt{3}i (- x_4 x_5 + x_4 x_6)+2 x_5^2 + 2 (\zeta_6 - 1) x_5 x_6 - 2 \zeta_6 x_6^2,$ \\
$-3 (-2 x_3 x_4 + x_4 x_5) + \sqrt{3}i (- x_3 x_5 + x_3 x_6 - x_4 x_6-x_5^2) + 2 \zeta_6 x_5 x_6 - x_6^2,$ \\
$-3 (2 x_3^2 + x_3 x_5 +x_4 x_6) +\sqrt{3}i (- x_3 x_6 + 3 x_4 x_5) +2 (\zeta_6 + 1) x_5 x_6 + 2 (\zeta_6 - 1) x_6^2$ \end{tabular}\\
\begin{tabular}{l}Maps:  \\
$(x_0,\ldots,x_6) \mapsto  (\zeta_3x_0,-\zeta_6x_2,-\zeta_3x_1-x_2,-x_3+\zeta_3x_4,\zeta_6x_3,\zeta_3 x_5,-\zeta_6 x_5+x_6)$, \\
$(x_0,\ldots,x_6) \mapsto  (-x_0,-x_2,x_1,x_4,-x-3,\zeta_3x_5-\zeta_6x_6,-\zeta_6x_5-\zeta_3x_6)$\end{tabular}\\
\mbox{} \\
Genus 7, Locus 10:  Group $(42,4)$, signature (2,6,21) cyclic trigonal\\
\begin{tabular}{l}Trigonal equation: $y^3 = x^8-x$\end{tabular}\\
\begin{tabular}{ll}Ideal: & $2\times 2$ minors of 
 $\left[ \begin{array}{rrrrr} x_0 & x_2 & x_3 & x_4 & x_5 \\ x_1 & x_3 & x_4 & x_5 & x_6 \end{array} \right]$, and \\
& $ x_0^3-x_6^2 x_2+x_2^2 x_3$ \\
& $x_0^2 x_1-x_6^2 x_3+x_2^2 x_4,$ \\
& $x_0 x_1^2-x_6^2 x_4+x_2^2 x_5,$ \\
& $x_1^3-x_6^2 x_5+x_2^2 x_6$ \end{tabular}\\
\begin{tabular}{l}Maps:\\
  $(x_0,\ldots,x_6) \mapsto  (\zeta_7 x_1,\zeta_7^6 x_0,-\zeta_7^4 x_6,-\zeta_7^2x_5,-x_4,-\zeta_7^5 x_3,-\zeta_7^3 x_2)$, \\
 $(x_0,\ldots,x_6) \mapsto  (\zeta_{21}^{-2} x_1, \zeta_{21}^{-5} x_0,-\zeta_{21}^{-1} x_6, -\zeta_{21}^{-4} x_5, (\zeta_3+1) x_4, -\zeta_{21}^{11} x_3, -\zeta_{21}^8 x_2 )$\end{tabular}\\
\mbox{} \\
Genus 7, Group 11:  Group $(32,11)$, signature (4,4,8). Has $g^{2}_{6}$ \\
\begin{tabular}{ll} Ideal: & $x_3 x_5+x_4 x_6,$ \\
& $ x_0^2+x_1 x_5+i x_2 x_6,$ \\
& $x_1 x_4+i x_2 x_3+x_5 x_6,$ \\
& $x_1 x_2+x_3 x_4,$ \\
& $x_1 x_6+\zeta_8^3 x_4 x_5,$ \\
& $x_2 x_5+\zeta_8 x_3 x_6,$ \\
& $x_1^2-i x_3^2-\zeta_8^3 x_5^2,$ \\
& $x_2^2+i x_4^2+\zeta_8^3 x_6^2,$ \\
& $-i x_2 x_4+\zeta_8^3 x_3^2,$ \\
& $x_1 x_3-\zeta_8^3 x_4^2$ \end{tabular}\\
\begin{tabular}{ll} Maps:& $(x_0,\ldots,x_6) \mapsto  (-x_0,-x_2,x_1,-ix_4,-ix_3,ix_6,ix_5)$, \\
& $(x_0,\ldots,x_6) \mapsto  (ix_0,x_1,ix_2,-x_3,-ix_4,-x_5,ix_6)$\end{tabular}\\
\mbox{} \\
Genus 7, Locus 12:  Group $(32,10)$, signature (4,4,8) \\
\begin{tabular}{ll}Ideal: &  $ x_1 x_6+\zeta_{16}^6 x_2 x_5+x_3 x_4, $ \\
& $x_1 x_2+x_5 x_6,$ \\
& $x_0^2+x_1 x_6-\zeta_{16}^6 x_2 x_5,$ \\
& $x_3 x_6-\zeta_{16}^4 x_4 x_5,$ \\
& $x_1^2-\zeta_{16}^7 x_4^2-\zeta_{16}^6 x_5^2,$ \\
& $x_2^2+\zeta_{16}^3 x_3^2-\zeta_{16}^{10} x_6^2,$ \\
& $-\zeta_{16}^2 x_2 x_6+(\zeta_{16}^{16}+\zeta_{16}^8) x_4^2-\zeta_{16}^7 x_5^2,$ \\
& $x_1 x_5+(-\zeta_{16}^{12}-\zeta_{16}^4) x_3^2-\zeta_{16}^{11} x_6^2,$ \\
& $x_1 x_3+\zeta_{16}^7 x_4 x_6,$ \\
& $x_2 x_4+\zeta_{16} x_3 x_5$ \end{tabular}\\
\begin{tabular}{ll} Maps: &$(x_0,\ldots,x_6) \mapsto  (-x_0,-x_2,x_1,-\zeta_{16}^2 x_4,-\zeta_{16}^6 x_3, -\zeta_{16}^6 x_6, -\zeta_{16}^2 x_5)$, \\
& $(x_0,\ldots,x_6) \mapsto (i x_0, -\zeta_{16}^6 x_2,-\zeta_{16}^2 x_1, -i x_4, i x_3, -i x_6, -i x_5)$\end{tabular}\\
Genus 7, Locus 13:  Group $(30,4)$, signature (2,15,30), hyperelliptic, $y^2=x^{15}-1$

\section*{References}
\bibliographystyle{amsplain}

\end{document}